
\magnification1200
\input amstex.tex
\documentstyle{amsppt}

\hsize=12.5cm
\vsize=18cm
\hoffset=1cm
\voffset=2cm

\footline={\hss{\vbox to 2cm{\vfil\hbox{\rm\folio}}}\hss}
\nopagenumbers
\def\DJ{\leavevmode\setbox0=\hbox{D}\kern0pt\rlap
{\kern.04em\raise.188\ht0\hbox{-}}D}

\def\txt#1{{\textstyle{#1}}}
\baselineskip=13pt
\def\hf{{\textstyle{1\over2}}}
\def\a{\alpha}\def\b{\beta}
\def\d{{\,\roman d}}
\def\e{\varepsilon}\def\E{{\roman e}}

\def\G{\Gamma}

\def\={\;=\;}

\def\zt{\zeta(\hf+it)}

\def\D{\Delta}
\def\no{\noindent}
 
\def\z{\zeta}

\def\hf{{\textstyle{1\over2}}}
\def\txt#1{{\textstyle{#1}}}

\def\le{\leqslant} \def\ge{\geqslant}
\font\tenmsb=msbm10
\font\sevenmsb=msbm7
\font\fivemsb=msbm5
\newfam\msbfam
\textfont\msbfam=\tenmsb
\scriptfont\msbfam=\sevenmsb
\scriptscriptfont\msbfam=\fivemsb

\font\ff=cmr8
\def\txt#1{{\textstyle{#1}}}
\baselineskip=13pt

\font\teneufm=eufm10
\font\seveneufm=eufm7
\font\fiveeufm=eufm5
\newfam\eufmfam
\textfont\eufmfam=\teneufm
\scriptfont\eufmfam=\seveneufm
\scriptscriptfont\eufmfam=\fiveeufm
\def\mathfrak#1{{\fam\eufmfam\relax#1}}

\font\tenmsb=msbm10
\font\sevenmsb=msbm7
\font\fivemsb=msbm5
\newfam\msbfam
     \textfont\msbfam=\tenmsb
      \scriptfont\msbfam=\sevenmsb
      \scriptscriptfont\msbfam=\fivemsb

  \def\rightheadline{{\hfil{\ff
  On zeta estimates in short intervals}\hfil\tenrm\folio}}

  \def\leftheadline{{\tenrm\folio\hfil{\ff
   Aleksandar Ivi\'c }\hfil}}
  \def\emptyheadline{\hfil}
  \headline{\ifnum\pageno=1 \emptyheadline\else
  \ifodd\pageno \rightheadline \else \leftheadline\fi\fi}

\topmatter
\title
ON SOME MEAN SQUARE ESTIMATES FOR THE ZETA-FUNCTION IN SHORT INTERVALS
\endtitle
\author   Aleksandar Ivi\'c  \endauthor

\bigskip
\address
Aleksandar Ivi\'c, Katedra Matematike RGF-a
Universiteta u Beogradu, \DJ u\v sina 7, 11000 Beograd, Serbia
\endaddress
\email {\tt
ivic\@rgf.bg.ac.rs,  aivic\@matf.bg.ac.rs} \endemail
\abstract
{Let $\D(x)$ denote the error term in the Dirichlet
divisor problem, and $E(T)$ the error term in the asymptotic
formula for the mean square of $|\zt|$. If
$E^*(t) = E(t) - 2\pi\D^*(t/2\pi)$ with $\D^*(x) =
 -\D(x)  + 2\D(2x) - \hf\D(4x)$ and we set
 $\int_0^T E^*(t)\d t = 3\pi T/4 + R(T)$, then we obtain
$$
\int_T^{T+H}(E^*(t))^2\d t \gg HT^{1/3}\log^3T
$$
and
$$
HT\log^3T \ll \int_T^{T+H}R^2(t)\d t \ll HT\log^3T,
$$
for $T^{2/3+\e}\le H \le T$. }
\endabstract
\endtopmatter

\document
\head
1. Introduction and statement of results
\endhead

This paper is the continuation of the author's works [6], [7], where the analogy
between the Riemann zeta-function $\z(s)$ and the divisor problem
was investigated.
As usual, let the error term in the classical Dirichlet divisor problem be
$$
\D(x) \;=\; \sum_{n\le x}d(n) - x(\log x + 2\gamma - 1),
\leqno(1.1)
$$
and
$$
E(T) \;=\;\int_0^T|\zt|^2\d t - T\left(\log\bigl({T\over2\pi}\bigr) + 2\gamma - 1
\right),\leqno(1.2)
$$
where $d(n)$ is the number of divisors of
$n$, $\z(s)$ is the Riemann zeta-function, and $ \gamma = -\G'(1) = 0.577215\ldots\,$
is Euler's constant. In view of F.V. Atkinson's classical explicit formula
for $E(T)$ (see [1] and [3, Chapter 15]) it was known long ago that
there are analogies between $\D(x)$ and $E(T)$. However, if one wants to stress
the analogy between $\z^2(s)$ and the divisor function, then
instead of the error-term function $\D(x)$ it
is more exact to work with the modified
function $\D^*(x)$ (see  M. Jutila [8], [9] and T. Meurman [10]), where
$$
\eqalign{
\D^*(x) :&= -\D(x)  + 2\D(2x) - \hf\D(4x)\cr&
= \hf\sum_{n\le4x}(-1)^nd(n) - x(\log x + 2\gamma - 1),\cr}
\leqno(1.3)
$$
since it turns out that $\D^*(x)$ is a better analogue of $E(T)$ than $\D(x)$.
Namely, M. Jutila (op. cit.) investigated both the
local and global behaviour of the difference
$$
E^*(t) \;:=\; E(t) - 2\pi\D^*\bigl({t\over2\pi}\bigr),
$$
and in particular in [9] he proved that
$$
\int_T^{T+H}(E^*(t))^2\d t \;\ll_\e\; HT^{1/3}\log^3T+ T^{1+\e}\quad(1\le H\le T).\leqno(1.4)
$$
Here and later $\e$ denotes positive constants which are arbitrarily
small, but are not necessarily the same ones at each occurrence,
while $a \ll_\e b$ (same as $a = O_\e(b))$ means that
the $\ll$--constant depends on $\e$. The significance of (1.4) is that, in view of
(see e.g., [3])
$$
\int_0^T(\D^*(t))^2\d t \sim AT^{3/2},\quad \int_0^T E^2(t)\d t \sim BT^{3/2}\quad(A,B >0,
T\to\infty),
$$
it transpires that $E^*(t)$ is in the mean square sense of a lower order of magnitude than
either $\D^*(t)$ or $E(t)$.

In [7] the author sharpened (1.4) (in the case when $H=T$) to the asymptotic formula
$$
\int_0^T (E^*(t))^2\d t \;=\; T^{4/3}P_3(\log T) + O_\e(T^{7/6+\e}),\leqno(1.7)
$$
where $P_3(y)$ is a polynomial of degree three in $y$ with
positive leading coefficient, and all the coefficients may be evaluated
explicitly.  This, in particular, shows that (1.4) may be complemented with the
lower bound
$$
\int_T^{T+H}(E^*(t))^2\d t \;\gg\; HT^{1/3}\log^3T\quad(T^{5/6+\e} \le H \le T).\leqno(1.8)
$$
It seems likely that the error term in (1.7) is $O_\e(T^{1+\e})$, but this seems difficult
to prove.

\medskip
In [6] the author investigated higher moments of $E^*(t)$, and e.g., in the second part of [6]
he proved that
$$
\int_0^T (E^*(t))^5\d t \;\ll_\e\; T^{2+\e};\leqno(1.9)
$$
but neither (1.4) nor (1.9) seem to imply each other.

\medskip
In part III of [6]  the error-term function $R(T)$ was  introduced by
 the relation
$$
\int_0^T E^*(t)\d t = {3\pi\over 4}T + R(T).\leqno(1.10)
$$
It was shown, by using an estimate for two-dimensional exponential sums, that
$$
R(T) = O_\e(T^{593/912+\e}), \quad {593\over912} = 0.6502129\ldots\,.
\leqno(1.11)
$$
It was  also proved that

\bigskip

$$
\int_0^TR^2(t)\d t = T^2p_3(\log T) + O_\e(T^{11/6+\e}),\leqno(1.12)
$$
where $p_3(y)$ is a cubic polynomial in $y$ with positive leading
coefficient, whose all coefficients may be explicitly evaluated,
and
$$
\int_0^TR^4(t)\d t \ll_\e T^{3+\e}.
\leqno(1.13)
$$

\medskip
The asymptotic formula (1.12) bears resemblance to (1.7), and it is proved
by a similar technique. The exponents in the error terms are, in both cases,
less than the exponent of $T$ in the main term by 1/6. From (1.7) one obtains that $E^*(T)
= \Omega(T^{1/6}(\log T)^{3/2})$, which shows that $E^*(T)$ cannot be
too small ($f = \Omega(g)$ means that $f = o(g)$ does not hold). Likewise, (1.12) yields
$$
R(T) \= \Omega\Bigl(T^{1/2}(\log T)^{3/2}\Bigr).\leqno(1.14)
$$
It seems plausible that the error term in (1.12) should be $O_\e(T^{5/3+\e})$,
while (1.14) leads one to suppose that
$$
R(T) = O_\e(T^{1/2+\e})\leqno(1.15)
$$
holds.

\bigskip
The aim of this paper to prove the following results.

\medskip
THEOREM 1. {\it For $T^{2/3+\e} \le H \le T$ we have}
$$
\int_T^{T+H}(E^*(t))^2\d t \;\gg\; HT^{1/3}\log^3T.\leqno(1.16)
$$

\medskip
\no
Note that (1.16) improves the range of $H$ for which (1.8) holds. \medskip
THEOREM 2. {\it For $T^{2/3+\e} \le H \le T$ we have}
$$
\int_T^{T+H}R^2(t)\d t \;\gg \; HT\log^3T,\leqno(1.17)
$$
{\it and, for $T^\e \le H \le T$,}
$$
\int_T^{T+H}R^2(t)\d t \;\ll_\e\; HT\log^3T + T^{5/3+\e}.\leqno(1.18)
$$

\medskip\no
The range for which (1.17) holds improves on the range for which  (1.12) holds.
\medskip
{\bf Corollary}. If $H = T^{2/3+\e}$, then every interval $[T, T+H]\;(T\ge T_0)$ contains
points $t_1, t_2$ such that, for some positive constants $A,B >0$,
$$
|E^*(t_1)| > At_1^{1/6}\log^{3/2}t_1,\quad |R(t_2)| > Bt_2^{1/2}\log^{3/2}t_2.\leqno(1.19)
$$
Note that this result follows from the asymptotic formulas (1.7) and (1.12), but in the
poorer range $T^{5/6+\e} \le H \le T$.
It would be interesting to find large positive and large negative values for which
the analogues of (1.19) hold. This was done in [4] for $E(t)$ and $\D(x)$, where it was shown
that there exist two positive constants $C,D$ such that, for $T\ge T_1$, every interval
$[T, T+C\sqrt{T}]$ contains points $t_3, t_4, t_5, t_6$ such that
$$
E(t_3) > D\sqrt{t_3},\; E(t_4) < -D\sqrt{t_4}, \;\D(t_5) > D\sqrt{t_5},
\; \D(t_6) <- D\sqrt{t_6}.\leqno(1.20)
$$
It would be interesting to obtain the analogue of (1.19) for large positive and negative values
of $E^*(t)$ and $R(t)$, like we have it in (1.20) for $E(t)$ and $\D(t)$, but this seems difficult.
\head
2. The necessary lemmas
\endhead

\bigskip
In this section we shall state the lemmas which are necessary
for the proof of our theorems.
 The first
two are Atkinson's classical explicit
formula for $E(t)$ (see e.g., [2] or [3])
and the Vorono{\"\i}-type formula for $\D^*(x)$, which is the analogue
of the classical truncated Vorono{\"\i} formula for $\D(x)$ (see [10]). The third
is an asymptotic formula involving $d^2(n)$.

\medskip
LEMMA 1. {\it Let $0 < A < A'$ be any two fixed constants
such that $AT < N < A'T$, and let $N' = N'(T) =
T/(2\pi) + N/2 - (N^2/4+ NT/(2\pi))^{1/2}$. Then }
$$
E(T) = \Sigma_1(T) + \Sigma_2(T) + O(\log^2T),\leqno(2.1)
$$
{\it where}
$$
\Sigma_1(T) = 2^{1/2}(T/(2\pi))^{1/4}\sum_{n\le N}(-1)^nd(n)n^{-3/4}
e(T,n)\cos(f(T,n)),\leqno(2.2)
$$
$$
\Sigma_2(T) = -2\sum_{n\le N'}\frac{d(n)}{n^{1/2}(\log T/(2\pi n))}
\cos\left(T\log \Bigl( {T\over2\pi n}\Bigr) - T + {1\over4}\pi\right),
\leqno(2.3)
$$
{\it with}
$$
\eqalign{\cr&
f(T,n) = 2T{\roman {arsinh}}\,\bigl(\sqrt{\pi n/(2T)}\,\bigr) + \sqrt{2\pi nT
+ \pi^2n^2} - {\txt{1\over4}}\pi\cr&
=  -\txt{1\over4}\pi + 2\sqrt{2\pi nT} +
\txt{1\over6}\sqrt{2\pi^3}n^{3/2}T^{-1/2} + a_5n^{5/2}T^{-3/2} +
a_7n^{7/2}T^{-5/2} + \ldots\,,\cr}\leqno(2.4)
$$
$$\eqalign{\cr
e(T,n) &= (1+\pi n/(2T))^{-1/4}{\Bigl\{(2T/\pi n)^{1/2}
{\roman {arsinh}}\,\Bigl(\sqrt{\pi n/(2T)}\,\Bigr)\Bigr\}}^{-1}\cr&
= 1 + O(n/T)\qquad(1 \le n < T),
\cr}\leqno(2.5)
$$
{\it and $\,{\roman{arsinh}}\,x = \log(x + \sqrt{1+x^2}\,).$}

\medskip
LEMMA 2 (see [2, Chapter 15]).{\it We have, for $1\ll N\ll x$},
$$
\D^*(x) = {1\over\pi\sqrt{2}}x^{1\over4}
\sum_{n\le N}(-1)^nd(n)n^{-{3\over4}}
\cos(4\pi\sqrt{nx} - {\txt{1\over4}}\pi) +
O_\e(x^{{1\over2}+\e}N^{-{1\over2}}).
\leqno(2.6)
$$
\medskip

LEMMA 3. {\it For $a > -\hf$ a constant we have}
$$
\sum_{n\le x}d^2(n)n^a = x^{a+1}P_3(\log x;a) + O_\e(x^{a+1/2+\e}),\leqno(2.7)
$$
{\it where $P_3(y;a)$ is a polynomial of degree three in $y$
whose coefficients depend on $a$,
and whose leading coefficient equals $1/(\pi^2(a+1))$. All the
coefficients of $P_3(y;a)$ may be explicitly evaluated}.

\medskip
This is a standard result, for a proof see e.g., Lemma 3 of [7].
\medskip

The next lemma  brings forth a formula
for $\int_0^T E(t)\d t$, which is closely related to
F.V.  Atkinson's classical explicit
formula for $E(T)$ (see [1] or e.g., Chapter 15 of [3] or Chapter 2 of [5]).
This is due to J.L. Hafner and the author [2] (see also Chapter 3 of [5]).

\bigskip
LEMMA 4. {\it We have}
$$
\eqalign{
\int_0^T E(t)\d t &= \pi T + {1\over2} \Bigl({2T\over\pi}\Bigr)^{3/4}
\sum_{n\le T}(-1)^nd(n)n^{-5/4}e_2(T,n)\sin f(T,n)\cr&
-2\sum_{n\le c_0T}d(n)n^{-1/2}\Bigl(\log {T\over2\pi n}\Bigr)^{-2}
\sin \left(T\log \Bigl( {T\over2\pi n}\Bigr) - T + {1\over4}\pi\right)
 \cr&+ O(T^{1/4}),\cr}\leqno(2.8)
$$
{\it where $\,c_0 = {1\over2\pi} + \hf - \sqrt{{1\over4}+{1\over2\pi}}\,,$
$\,{\roman{ar\,sinh}}\,x = \log(x + \sqrt{1+x^2}\,),$ and
for $\,1\le n \ll T$, }
$$
\eqalign{
e_2(T,n) &= \left(1+{\pi n\over T}\right)^{-1/4}\left\{\left({2T\over\pi n}\right)^{1/2}
{\roman {ar\,sinh}}\left({\pi n\over2T}\right)^{1/2}\right\}^{-1/2}\cr&
= 1 + b_1{n\over T} + b_2\Bigl({n\over T}\Bigr)^2 + \ldots\,,\cr
f(T,n) &= 2T{\roman {ar\,sinh}}\,\bigl(\sqrt{\pi n/(2T)}\,\bigr) + \sqrt{2\pi nT
+ \pi^2n^2} - {\txt{1\over4}}\pi\cr&
=  -\txt{1\over4}\pi + 2\sqrt{2\pi nT} +
a_3n^{3/2}T^{-1/2} + a_5n^{5/2}T^{-3/2} +
a_7n^{7/2}T^{-5/2} + \ldots\,.\cr}\leqno(2.9)
$$

\bigskip
We also need a formula for the integral of $\D^*(x)$. From a classical
result of G.F.  Vorono{\"\i}  [10] (this also easily follows from pp. 90-91
of [3]) we have
$$
\int_0^X\D(x)\d x = {X\over4} + {X^{3/4}\over2\sqrt{2}\pi^2}
\sum_{n=1}^\infty d(n)n^{-5/4}\sin(4\pi\sqrt{nX}- {\txt{1\over4}}\pi)
+ O(X^{1/4}).
$$
To relate the above integral to the one of $\D^*(x)$ we proceed as
on pp. 472-473 of [3], using (1.3) and (1.10). In this way we are
led to

\bigskip
LEMMA 5. {\it We have}
$$
\eqalign{
\int_0^T\D^*(t)\d t &= -{T\over8} +
 {T^{3/4}\over2\sqrt{2}\pi^2} \sum_{n\le T^{2}}
(-1)^nd(n)n^{-5/4}\sin(4\pi\sqrt{nT}- {\txt{1\over4}}\pi)\cr&
+ O(T^{1/4}).
\cr}\leqno(2.10)
$$

\bigskip

\head
3. Proof of Theorem 1
\endhead
\bigskip
We use Lemma 1 and Lemma 2
with $N=T$
to deduce that, for $T\le t\le T+H, T^{2/3+\e}\ll H\le T$,
$$
E^*(t) := S_1(t) + S_2(t) +S_3(t), \leqno(3.1)
$$
where
$$
\eqalign{
&S_1(t) :=\sqrt{2}\left({t\over2\pi}\right)^{1/4}
\sum_{n\le T}(-1)^nd(n)n^{-3/4}\Biggl\{e(t,n)\cos f(t,n)\cr&
\qquad\;- \cos (\sqrt{8\pi nt} - {\txt{1\over4}}\pi)\Biggr\}\cr&
S_2(t) := -2\sum_{n\le N'}d(n)n^{-1/2}\left(\log {t\over2\pi n}\right)^{-1}
\cos\left(t\log {t\over2\pi n}-t+ {\pi\over 4}\right),\cr&
S_3(t) := O_\e(T^\e),\cr}\leqno(3.2)
$$
and $N' = t/(2\pi) + t/2 - \sqrt{T^2/4 + tT/(2\pi)}$.
We have, similarly as in [7],
$$
\int_T^{T+H}\Bigl\{S_2^2(t) + S_3^2(t)\Bigr\}\d t \;\ll_\e\;HT^{\e},\leqno(3.3)
$$
since $S_2(t)$ is in fact quite analogous to the sum representing $\z^2(\hf+it)$.
Therefore
$$
\eqalign{&\int_T^{T+H}(E^*(t))^2\d t\cr&
= \int_T^{T+H}\Bigl\{S_1^2(t) + S_2^2(t) + S_3^2(t) + 2S_1(t)S_2(t)+2S_1(t)S_3(t)+2S_2(t)S_3(t)\Bigr\}\d t
\cr&
= \int_T^{T+H}S_1^2(t)\d t +  O_\e(HT^{1/6+\e}).\cr}\leqno(3.4)
$$
Here we used (3.3), (1.4) and the Cauchy-Schwarz inequality for integrals to deduce that
$$
\eqalign{
\int_T^{T+H}S_1(t)S_2(t)\d t &\ll \left\{\int_T^{T+H}S_1^2(t)\d t\int_T^{T+H}S^2_2(t)\d t\right\}^{1/2}\cr&
\ll_\e \left(HT^\e HT^{1/3}\log^3T\right)^{1/2} \ll_\e HT^{1/6+\e}.\cr}
$$
Now we write
$$
\eqalign{
S_1(t) &= S_4(t) + S_5(t),\cr
S_4(t)&:= \sqrt{2}\left({t\over2\pi}\right)^{1/4}\sum_{n\le T^{1/2-\e}}(-1)^nd(n)n^{-3/4}
\Biggl\{e(t,n)\cos f(t,n)- \cr&-\cos (\sqrt{8\pi nt} - {\txt{1\over4}}\pi)\Biggr\}\cr
S_5(t)&:= \sqrt{2}\left({t\over2\pi}\right)^{1/4}\sum_{T^{1/2-\e}<n\le T}(-1)^nd(n)n^{-3/4}
\Biggl\{e(t,n)\cos f(t,n)- \cr&-\cos (\sqrt{8\pi nt} - {\txt{1\over4}}\pi)\Biggr\}.\cr}
$$
We obtain, following the proof of (1.4),
$$
\int_T^{T+H}S_1^2(t)\d t = \int_T^{T+H}\left\{S_4^2(t) + S_5^2(t) + 2S_4(t)S_5(t)\right\}\d t.
$$
In view of (1.4) we have
$$
\int_T^{T+H}S_4^2(t)\d t \;\ll_\e\;HT^{1/3+\e}\qquad(T^{2/3+\e}\le H\le T).\leqno(3.5)
$$
To estimate the mean square of $S_5(t)$, we split the sum into subsums with the range of
summation $K<n\le K'\le 2K, T^{1/2-\e}\le K \le T$.
Note that the mean square bound ($c\ne0$)
$$
\eqalign{&
\int_T^{T+H}\Bigl|\sum_{K< k\le K'\le2K}(-1)^kd(k){\roman e}^{\sqrt{ckt}i}\Bigr|^2\d t
\cr&
= H\sum_{K< k\le2K}d^2(k) + \sum_{K< m\ne n\le2K}(-1)^{m+n}d(m)d(n)
\int_T^{T+H}{\roman e}^{\sqrt{ct}(\sqrt{m}-\sqrt{n})i}\d t\cr&
\ll HK\log^3T + \sqrt{T}\sum_{K< m\ne n\le2K}{d(m)d(n)\over|\sqrt{m}-\sqrt{n}|}\cr&
\ll_\e HK\log^3T + T^{1/2+\e}\sum_{K< m\ne n\le2K}{K^{1/2}\over|m-n|}\cr&
\ll_\e T^\e(HK + T^{1/2}K^{3/2}) \cr}\leqno(3.6)
$$
holds for $1\ll K \ll T^C\;(C>0)$, where we used the standard first derivative test
for exponential sums
(see Lemma 2.1 of [3]) and Lemma 3. The same bound also holds if in the exponential we have
$f(t,k)$ (cf. (2.4)) instead of $\sqrt{ctk}$, as shown e.g.,
in the derivation of the mean square
formula for $E(t)$ in Chapter 15 of [3]. Using (3.6) it follows that
$$
\int_T^{T+H}S_5^2(t)\d t \;\ll_\e\; T^{1/2+\e}(HT^{-1/4}+T^{1/2}) = HT^{1/4+\e}+T^{1+\e}\leqno(3.7)
$$
holds for $T^{2/3+\e}\le H\le T$. Consequently using (3.5), (3.7) and the Cauchy-Schwarz inequality
we obtain
$$
\int_T^{T+H}S_4(t)S_5(t)\d t\;\ll_\e\;T^\e(HT^{7/24}+ H^{1/2}T^{2/3}).
$$
Therefore, for $T^{2/3+\e}\le H\le T$, we have shown that
$$
\int_T^{T+H}S_1^2(t)\d t = \int_T^{T+H}S_4^2(t)\d t +
 O_\e\Bigr(T^\e(HT^{7/24}+ H^{1/2}T^{2/3})\Bigl).
\leqno(3.8)
$$
The integral on the right-hand side of (3.8) is equal to
$$
\eqalign{&
\sqrt{{2\over\pi}}\sum_{n\le T^{1/2-\e}}d^2(n)n^{-3/2}\int_T^{T+H}t^{1/2}\left(e(t,n)\cos f(t,n)
- \cos(\sqrt{8\pi nt}-\pi/4)\right)^2\d t\cr&
+ \sqrt{{2\over\pi}}\sum_{1\le m\ne n\le T^{1/2-\e}}(-1)^{m+n}d(m)d(n)(mn)^{-3/4}\times\cr&
\times\int_T^{T+H}t^{1/2}\Bigl(e(t,m)\cos f(t,m)-\ldots\Bigr)\Bigl(e(t,n)\cos f(t,n)-\ldots)\Bigr)\d t.\cr}
$$
In this expression we first replace the factors $e(t,m)$ and $e(t,n)$ by 1, and it is seen that the
total error made in this process is $O_\e(HT^{1/4+\e})$, since $e(t,n) = 1 + O(n/t)$ and
$m, n\le T^{1/2-\e}$. Consider now the sum over $m\ne n$. If both $m$ and $n$ are $\le T^{1/3-\e}$,
then observe that Taylor's formula gives
$$
\eqalign{
& \sin f(t,m) - \sin(2\sqrt{2\pi mt} - \pi/4) = \sum_{k=1}^\infty
{(y-y_0)^k\over k!}\sin(y_0 + \hf k\pi)\cr&
y = f(t,m),\; y_0 = 2\sqrt{2\pi mt} - \pi/4, \; y-y_0= d_3m^{3/2}t^{-1/2}
+ d_5m^{5/2}t^{-3/2} + \ldots\,,\cr}\leqno(3.9)
$$
and similarly for $\sin f(t,n)$. Therefore the total contribution of these terms, by the first derivative test,
will be
$$
\eqalign{&
\ll T\sum_{m\ne n\le T^{1/3-\e}}d(m)d(n)(mn)^{3/4}T^{-1}\cdot\frac{1}{|\sqrt{m}-\sqrt{n}|}\cr&
\ll_\e T^\e\sum_{m\ne n\le T^{1/3-\e}}\frac{(mn)^{3/4}(\sqrt{m}+\sqrt{n})}{|m-n|} \ll_\e T^\e
\sum_{m\le T^{1/3-\e}}m^2 \ll_\e T^{1+\e}.\cr}
$$
If, say, $m\le T^{1/3-\e}, T^{1/3-\e}<n\le T^{1/2-\e}$, then the contribution is a multiple of
$$
\eqalign{&
\sum_{m\le T^{1/3-\e}}(-1)^md(m)m^{3/4}\sum_{T^{1/3-\e}<n\le T^{1/2-\e}}(-1)^nd(n)n^{-3/4}\times\cr&
\times \int_T^{T+H}\E^{\pm i\sqrt{8\pi mt}}\left(\E^{\pm if(t,n)}-\E^{\pm i\sqrt{8\pi nt}}\right)\d t.\cr}
$$
The contribution of terms with two square roots in the exponential is, by the first derivative test,
$$
\eqalign{&
\ll\sum_{m\le T^{1/3-\e}}d(m)m^{3/4}\sum_{T^{1/3-\e}<n\le T^{1/2-\e}}d(n)n^{-3/4}\frac{T^{1/2}}{|\sqrt{m}-\sqrt{n}|}
\cr&
\ll_\e T^{3/4+\e}\sum_{T^{1/3-\e}<n\le T^{1/2-\e}}d(n)n^{-1/2} \ll_\e T^{1+\e}.\cr}\leqno(3.10)
$$
The remaining  case of interest is when we have the exponential factor
$$
\exp\Bigl(iF(t,m,n)\Bigr),\quad F(t,m,n) := \sqrt{8\pi mt} - f(t,n),
$$
when
$$
\eqalign{
{\d\over\d t}F(t,m,n) &= \sqrt{2\pi m\over t} - 2{\roman{ar\,sinh}}\sqrt{\pi n\over 2t}\cr&
= \sqrt{2\pi\over t}(\sqrt{m} - \sqrt{n}\,) + a_3n^{3/2}t^{-3/2} + a_5n^{5/2}t^{-5/2}
+ \ldots\,.\cr}
$$
But as, for $m\ne n$ and $m,n \le T^{1/2-\e}$,
$$
\frac{|\sqrt{m} - \sqrt{n}|}{\sqrt{t}} = \frac{|m-n|}{|\sqrt{m}+\sqrt{n}|\sqrt{t}} \ge
\frac{1}{(\sqrt{m}+\sqrt{n})\sqrt{t}} \gg \frac{t^\e n^{3/2}}{t\sqrt{t}},
$$
then again by the first derivative test we obtain a contribution which is, similarly to (3.8), $\ll_\e T^{1+\e}$.
Finally, the same argument shows that the contribution, when $T^{1/3-\e} < m\ne n \le T^{1/2-\e}$, is
$$
\eqalign{&
\ll T^{1/2}\sum_{T^{1/3-\e}<m\le T^{1/2-\e}}d(m)m^{-3/4}\sum_{T^{1/3-\e}<n \le T^{1/2-\e}, n\ne m}d(n)n^{-3/4}
\frac{\sqrt{T}}{|\sqrt{m}-\sqrt{n}|}\cr&
\ll_\e T^{1+\e}\sum_{T^{1/3-\e} < m\ne n \le T^{1/2-\e}}\frac{(mn)^{-3/4}(\sqrt{m}+\sqrt{n})}{|m-n|}
\cr&
\ll_\e T^{1+\e}\sum_{T^{1/3-\e} < m \le T^{1/2-\e}}m^{-1}\log T \ll_\e T^{1+\e}.\cr}
$$
From (3.8) and the preceding  estimates it follows that, for $T^{2/3+\e} \le H \le T$,
$$
\eqalign{&
\int_T^{T+H}(E^*(t))^2\d t =  O_\e(HT^{7/24+\e} + H^{1/2}T^{2/3+\e} + T^{1+\e}) +\cr&
+\,\sqrt{2\over\pi}\sum_{n\le T^{1/2-\e}}\frac{d^2(n)}{n^{3/2}}
\int_T^{T+H}t^{1/2}\Bigl\{\cos f(t,n)-\cos(\sqrt{8\pi nt}-\pi/4)\Bigr\}^2\d t.\cr}
\leqno(3.11)
$$
Since the integrand on the right-hand side of (3.11) is non-negative, it is not difficult to deduce
(1.16) of Theorem 1 from (3.11). To manage the cosines in (3.11) we use the elementary
identity
$$
(\cos\a-\cos\b)^2 = 1 - \cos(\a-\b)+\{\hf\cos2\a +\hf\cos2\b - \cos(\a+\b)\}
$$
with $\a = f(t,n),\b = \sqrt{8\pi nt} - {\txt{1\over4}}\pi$.
By the first derivative test it is seen that the terms coming from curly braces
contribute $\ll T$ to (3.11). Furthermore, in view of
$$
1-\cos \gamma = 2\sin^2(\hf \gamma), \;|\sin x| \ge {2\over\pi}|x|\quad(|x| \le \pi/2),
$$
it is seen that the sum on the right-hand side of (3.11) is
$$
\eqalign{&
\sqrt{8\over\pi}\sum_{n\le T^{1/2-\e}}d^2(n)n^{-3/2}\int_T^{T+H}t^{1/2}
\sin^2\Bigl(a_3n^{3/2}t^{-1/2}+a_5n^{5/2}t^{-3/2}+\cdots\Bigr)\d t + \cr& +O_\e(T^{1+\e})\cr&
\gg
\sum_{n\ll T^{1/3}}d^2(n)n^{-3/2}\int_T^{T+H}t^{1/2}n^3t^{-1}\d t + O_\e(T^{1+\e})\cr&
\ge C_1T^{-1/2}H\sum_{n\ll T^{1/3}}d^2(n)n^{3/2}+ O_\e(T^{1+\e})\cr&
\ge C_2T^{-1/2}HT^{5/6}\log^3T  + O_\e(T^{1+\e}) \ge \hf C_2HT^{1/3}\log^3T\cr}
$$
for $T^{2/3+\e} \le H \le T$. Since all the $O$-terms in (3.11) are $o(HT^{1/3}\log^3T)$
in this range, it means that we have proved (1.16) of Theorem 1.

\medskip
\head
4. Proof of Theorem 2
\endhead
\medskip
Combining Lemma 4 and Lemma 5 we obtain, with $c_0$ as in (2.8),
$$
\eqalign{
&R(T) =
{1\over2}\left({2T\over\pi}\right)^{3/4}
\sum_{T<n\le T^{2}}(-1)^{n+1}d(n)n^{-5/4}\sin(2\pi\sqrt{2nT}- {\txt{1\over4}}\pi)\cr&
+ {1\over2}\left({2T\over\pi}\right)^{3/4} \sum_{n\le T}
(-1)^nd(n)n^{-5/4}\left\{e_2(T,n)\sin f(T,n)-\sin(2\pi\sqrt{2nT}-
{\txt{1\over4}}\pi)\right\}\cr&
-2\sum_{n\le c_0T}d(n)n^{-1/2}\Bigl(\log {T\over2\pi n}\Bigr)^{-2}
\sin \left(T\log \Bigl( {T\over2\pi n}\Bigr) - T + {\txt{1\over4}}\pi\right)
+ O(T^{1/4}).
\cr}\leqno(4.1)
$$
This gives, since the estimation of $\sum_{n\le N'}\cdots$ is similar (see e.g., [3])
to the estimation of 
$$
\z^2(\hf+it) \;=\;O_\e(t^{1/3}),
$$
$$
\eqalign{
&R(T) = O(T^{1/2}\log T)\, +\cr&
+ {1\over2} \Bigl({2T\over\pi}\Bigr)^{3/4}
\sum_{n\le T}(-1)^nd(n)n^{-5/4}
\left\{e_2(T,n)\sin f(T,n)-\sin(2\sqrt{2\pi nT}-\pi/4)
\right\}.\cr}
\leqno(4.2)
$$
We further simplify (4.2) by estimating trivially the portion of the sum for which $n > T^{1/2-\e}$ and
then using $e_2(T,n) = 1 + O(n/T)$. We obtain
$$
\eqalign{
&R(T) = O_\e(T^{1/2+\e})\, +\cr&
+ {1\over2} \Bigl({2T\over\pi}\Bigr)^{3/4}
\sum_{n\le T^{1/2-\e}}(-1)^nd(n)n^{-5/4}
\left\{\sin f(T,n)-\sin(2\sqrt{2\pi nT}-\pi/4)\right\}\cr&
= {1\over2} \Bigl({2T\over\pi}\Bigr)^{3/4}(s_1(T) + s_2(T))
+ O_\e(T^{1/2+\e}),\cr}
\leqno(4.3)
$$
say, where in $s_1$ summation is over $n \le T^{1/3-\e}$, and in $s_2$ summation is over
$n$ such that $T^{1/3-\e} < n\le T^{1/2-\e}$.

\medskip
Now we replace $T$ by $t$ and suppose that $T \le t\le T+H, T^{2/3+\e} \le H \le T$.
We prove first (1.18) of Theorem 2. In $s_1(t)$ we use (3.9), and in $s_2(t)$ we consider
separately the contributions coming from $\sin f(t,n)$ and $\sin(\sqrt{8\pi nt}-\pi/4)$. In both
cases we use (3.6), since it was mentioned that the argument also works for $ f(t,n)$ in the exponential.
Thus we are led to the estimation of the integrals ($K\ll T^{1/3-\e}$)
$$
\eqalign{&
I_1(K) := \int_T^{T+H}T^{1/2}\Bigl|\sum_{K<n\le K'\le2K}(-1)^nd(n)n^{1/4}\E^{i\sqrt{8\pi nt}}\Bigr|^2\d t\cr&
\ll_\e \max_{K\ll T^{1/3-\e}}T^{1/2+\e}(HK^{3/2}+T^{1/2}K^2) \ll_\e T^{1+\e}H + T^{5/3+\e}\cr}
$$
and ($T^{1/3-\e} \ll K \ll T^{1/2-\e}$)
$$
\eqalign{&
I_2(K) := \int_T^{T+H}T^{3/2}\Bigl|\sum_{K<n\le K'\le2K}(-1)^nd(n)n^{-5/4}\E^{i\sqrt{8\pi nt}}\Bigr|^2\d t\cr&
\ll_\e \max_{T^{1/3-\e}\ll K\ll T^{1/2-\e}}T^{3/2+\e}(HK^{-3/2}+T^{1/2}K^{-1}) \ll_\e T^{1+\e}H + T^{5/3+\e},\cr}
$$
while the integral
$$
I_3(K) := \int_T^{T+H}T^{3/2}\Bigl|\sum_{K<n\le K'\le2K}(-1)^nd(n)n^{-5/4}\E^{if(t,n)}\Bigr|^2\d t
$$
is estimated analogously as $I_2(K)$.
Since
$$
\eqalign{
\int_T^{T+H}R^2(t) \d t &\ll_\e HT^{1+\e} + \log T\max_{K\ll T^{1/3-\e}}I_1(K)\cr&
+ \log T \max_{T^{1/3-\e}\ll K\ll T^{1/2-\e}}\Bigl(I_2(K) + I_3(K)\Bigr),\cr}
$$
the bound (1.18) follows.

\medskip
The proof of (1.17) is carried out by using (4.3) and (1.18), 
and is  analogous to the proof of Theorem 1, only somewhat less involved.
The sum corresponding to $S_4(t)$ in the proof of Theorem 1 (cf. (3.5))
is the main term on the right-hand side of (4.3).
There is no need to repeat the details.

\medskip


\topglue1cm
\bigskip
\Refs
\bigskip

\item{[1]} F.V. Atkinson, The mean value of the Riemann zeta-function,
Acta Math. {\bf81}(1949), 353-376.

\item{[2]} J.L. Hafner and A. Ivi\'c, On the mean square of the Riemann
zeta-function on the critical line, J. Number Theory {\bf31}(1989), 151-191.

\item{[3]} A. Ivi\'c, The Riemann zeta-function, John Wiley \&
Sons, New York, 1985 (2nd ed. Dover, Mineola, New York, 2003).

\item{[4]} A. Ivi\'c, Large values of certain number-theoretic error terms,
Acta  Arithm.  {\bf56}(1990), 135-159.

\item{[5]} A. Ivi\'c, The mean values of the Riemann zeta-function,
LNs {\bf 82}, Tata Inst. of Fundamental Research, Bombay (distr. by
Springer Verlag, Berlin etc.), 1991.

\item{[6]} A. Ivi\'c, On the Riemann zeta-function and the divisor problem,
Central European J. Math. {\bf(2)(4)} (2004), 1-15,  II, ibid.
{\bf(3)(2)} (2005), 203-214,  III, Annales Univ.
Sci. Budapest, Sect. Comp. {\bf29}(2008), 3-23,
and IV, Uniform Distribution Theory {\bf1}(2006), 125-135.

\item{[7]} A. Ivi\'c, On the mean square of the zeta-function and
the divisor problem, Annales  Acad. Scien. Fennicae Mathematica {\bf23}(2007), 1-9.

\item{[8]} M. Jutila, Riemann's zeta-function and the divisor problem,
Arkiv Mat. {\bf21}(1983), 75-96 and II, ibid. {\bf31}(1993), 61-70.

\item{[9]} M. Jutila, On a formula of Atkinson, Topics in classical number theory, Colloq. Budapest 1981, 
Vol. I, Colloq. Math. Soc. J\'anos Bolyai {\bf34}(1984), 807-823. 

\item{[10]} G.F. Vorono{\"\i}, Sur une fonction transcendante et ses applications
\`a la sommation de quelques s\'eries, Ann. \'Ecole Normale (3){\bf21}(1904),
2-7-267 and ibid. 459-533.

\endRefs
\vskip1cm

\enddocument

\bye